\documentclass[12pt]{IEEEtran}  

\usepackage{amsmath,amsthm,amsfonts,amssymb,color}
\usepackage{graphicx}
\usepackage[mathscr]{eucal}
\usepackage{float}
\usepackage{epstopdf}
\usepackage{fancyhdr}
\usepackage{latexsym}
\newtheorem{thm}{Theorem} 

\newtheorem{prob}[thm]{Problem} 
\newtheorem{lemma}{Lemma}

\newtheorem{rem}{Remark}

\newcommand{\beq}{\begin{equation}}
\newcommand{\eeq}{\end{equation}}

\newcommand{\nd}{\noindent}
\newcommand{\q}{\quad}
\newcommand{\qq}{\qquad}

\DeclareMathOperator{\tr}{tr}

\DeclareMathOperator{\diag}{diag}
\DeclareMathOperator{\rank} {rank}

\DeclareMathOperator{\im} {Im}

\newcommand{\Frac}[2]{{\displaystyle\frac{#1}{#2}}}

\newcommand{\bmat}{\left[ \begin{matrix}}
\newcommand{\emat}{\end{matrix} \right]}

\DeclareMathOperator{\E}{{\mathbb E}}

\newcommand{\Rbb}{\mathbb R}

\newcommand{\Zbb}{\mathbb Z}

\newcommand{\yb}{\mathbf y}
\newcommand{\zb}{\mathbf  z}

\newcommand{\Hb}{\mathbf H}

\def\xib{\boldsymbol{\xi}}

\def\Sigmab{\boldsymbol{\Sigma}}

\title{Approximation of stationary processes and Toeplitz Spectra}

\author{Giorgio Picci and Bin Zhu 
\thanks{Giorgio Picci is with the Department of Information Engineering, University of Padova, Padova, Italy {\tt\small picci@dei.unipd.it}}%
\thanks{$^{2}$ Bin Zhu is with the  School of Intelligent Systems Engineering, Sun Yat-sen University, Waihuan East Road 132, 510006 Guangzhou, China, e-mail: {\tt\small zhub26@mail.sysu.edu.cn}}%
 }

\begin{document}

\maketitle
\thispagestyle{empty}

\begin{abstract}

We study the approximation of stationary processes by a simple class of purely deterministic signals. This has an analytic counterpart in the approximation of symmetric positive definite Toeplitz matrices by  submatrices  of finite rank. We  propose a notion of distance between them and prove a weak sense approximation result. 
\end{abstract} 


\section{Introduction and Problem Statement}\label{Sec: ProbStat}
In our past investigations \cite{arkiv} we have discovered   that  the {\em a posteriori} covariance operator of some random  harmonic oscillation signals  corrupted by white noise  has remarkable properties very similar to those that have been   uncovered  in the 60's and 70's by D.~Slepian and coworkers in a famous series of papers concerning the energy concentration problems of time and band limited signals \cite{Slepian-P-61,Slepian-P-61bis,Slepian-78,Landau-W-80}. The key  property of the covariance  operator in question is that its eigenvalues decay extremely fast to zero for indices greater than an   a priori computable number (the so-called {\em Slepian frequency}).

From the sharp decay to zero of the eigenvalues it  follows that the a posteriori probabilistic description of the signal is essentially finite dimensional and it seems that in simulated experiments the observations can be well approximated by a {\em purely deterministic process} \cite{LPBook} corrupted by white noise.  Since purely deterministic processes of finite rank can be described by linear finite dimensional state-space models   \cite[p. 277]{LPBook}, the estimation can be carried on by rather straightforward {\em subspace methods}.

This although the precise meaning of this approximation has so far not been  understood. Scope of this paper is to propose a possible metric for this approximation and to prove  convergence althogh only in a weak sense.

Consider the infinite covariance  matrix\footnote{In this paper  bold symbols like  $\Sigmab$ are reserved for   {\bf infinte} arrays say stochastic processes or  covariance matrices.}
$
\Sigmab= \left [ \sigma(t-s)\right]_{t,s \in \Zbb}
$ of a  scalar stationary zero-mean  process having covariance function
\begin{equation}\label{Model}
\sigma(\tau) = \E\yb (t+\tau)\yb (t) \, \qquad \tau \in \Zbb
\end{equation}
which we shall assume   admits a Fourier transform
$$
\varphi(e^{j\theta}) = \sum_{\tau=-\infty}^{+\infty}\, e^{-j\theta\, \tau} \sigma(\tau)
$$
which is a piecewise smooth function of $\theta$. For example $\varphi$ will be continuous and bounded when $\sigma$ belongs to $\ell^1$. The function $\varphi$ is called the  {\em symbol} of $\Sigmab$. The process $\yb$ is not necessarily purely non deterministic; however  we assume that $\varphi(e^{j\theta}) $ is piecewise continuous and bounded, as for example is a rectangular function. Then  $\Sigmab$ induces  a bounded operator in $\ell^2$ \cite{Akhiezer-G-61}.

Let now 
$$\yb_N(t)= \bmat \yb(t)&\yb(t+1)&\ldots&\yb(t+N-1)\emat^{\top}
$$ and consider the $N$-truncation of the matrix $\Sigmab$, defined as
$$
\Sigma_N := \E\yb_N(t)\yb_N(t)^{\top}
$$
which for each $N$ has a positive point spectrum, say
$$
\frak{S}_N:= \{	\lambda_{N,1}, \ldots, \lambda_{N,N}\}
$$
where the eigenvalues are ordered in decreasing magnitude. We shall assume that all $\Sigma_N$ are non singular so that the eigenvalues must be strictly positive for all $N$. By Weyl's theorem \cite[p. 203]{Stewart-S-90}, see also   \cite[Fact M]{forni_lippi_2001},     the k--th eigenvalue of $\Sigmab_N$ is a non decreasing function of $N$ and hence has a limit, $\lambda_k(\Sigmab)$, which may possibly be equal to  $+\infty$. Each such limit is called an {\em  eigenvalue of $\Sigmab$}.
These limits however  are in general not true eigenvalues, as it is well-known that $\Sigmab$ may not have eigenvalues. For example, a bounded symmetric Toeplitz operator matrix has a purely continuous spectrum \cite{Hartman-W-54}. 

Assume now that all eigenvalues of $\Sigmab$ are finite (this is equivalent to $\Sigmab$ being a bounded linear operator in $\ell^2$ see \cite{Bottegal-P-14}) and let us keep in the formal spectral decomposition of $\Sigmab$ the  $n$ largest eigenvalues setting all the others to zero. In this way we form an "approximation" of $\Sigmab$ of finite rank $n$ which we shall denote $\Sigmab^n$. Then  the (infinite) matrix $\Sigmab^n$ is the covariance   of a purely deterministic process of rank $n$ \cite{LPBook} whose spectral density, say $\varphi_n(e^{j\theta}) $, is the sum of $n$ Dirac pulses. We would like to 	know in what sense, if any,  $\Sigmab^n$ could be considerd an approximation of $\Sigmab$ or equivalently, $\varphi_n(e^{j\theta})$ could  be considered an approximation of the symbol $\varphi(e^{j\theta})$. Obviously this last     fact  could only happen in a weak sense, say for arbitrary test functions $\psi(e^{j\theta})$ continuous on  the unit circle one should have
\begin{equation}\label{Wconv}
\int  \psi(e^{j\theta})\varphi_n(e^{j\theta}) \psi(e^{j\theta})^*\to  \int  \psi(e^{j\theta})\varphi(e^{j\theta}) \psi(e^{j\theta})^*
\end{equation}
as $n\to \infty$.

An equivalent question can be posed in terms of $L^2$ approximation of a  stationary  process $\yb$ of covariance $\Sigmab$ (which could in particular be p.n.d) by a purely deterministic one, say $\hat\yb_n$ having covariance $\Sigmab_n$ of rank $n$. As we shall see this problem  should   also be naturally  formulated in a weak sense.

\section{Approximation of random vectors}
To begin with, suppose we want to  approximate an $N$-dimensional  zero-mean random vector $y$ of covariance matrix   $\Sigma \in \Rbb^{N\times N}$ by another $N$-dimensional vector say $\hat{y}$ having covariance $\Sigma_n$ of rank $n<N$. To make the problem well-posed we shall require that the approximation  $\hat{y}$ generates a linear  subspace  of  $\Hb({y})$, which will then have to be $n$-dimensional. This means that we can represent $\hat{y}$ as a linear function of $y$, say
$$
\hat{y}:= M {y}
$$
where $M\in \Rbb^{N\times N}$ has rank $n$.

Motivated from the above,    let us consider the following approximation problem:
\begin{prob}\label{ApproxProb}
Find a   matrix  $M \in {\mathbb R}^{N\times N}$ of rank $n$, solving the following minimum problem
\begin{equation}  \label{Approx}
\min_{\rank (\,M\,)\, =\, n}\, \E\{\|{y} -M\, {y}\|^2\}\,.
\end{equation}
\end{prob}
Note that an equivalent geometric formulation is  to look for an optimal $n$-dimensional subspace of  $\Hb({y})$ onto which ${y}$ should be projected in order to minimize the approximation error variance. Let us  stress that this is quite different from the usual least squares approximation problem which amounts to projecting onto a {\em given subspace}.

 As usual, minimizing the  square distance in \eqref{Approx} requires that the approximation $ M {y}$ should be uncorrelated with the approximation error; namely
\begin{equation}  \label{OrthM}
{y} - M{y} \,\perp\, M{y} 
\end{equation}
which is equivalent to
$$
M\Sigma -M \Sigma M^{\top} = 0\,.
$$
Introducing a  square root $ \Sigma^{1/2}$ of $\Sigma$ and defining $\hat M := \Sigma^{-1/2} M
\Sigma^{1/2}$, this condition is seen to be equivalent to 
$$
\hat M = \hat M \,\hat M^{\top}
$$
which just says that $ \hat M$ must be {\it  symmetric and  idempotent} (i.e. $\hat M = \hat M^2$), in other words an  {\it
orthogonal projection} from ${\mathbb R}^{N}$ onto some  $n$-dimensional subspace. Hence  $M$ must have the following structure
\begin{equation} \label{structM}
M\,=\, \Sigma^{1/2} \,\Pi\,\Sigma^{-1/2},\quad\quad \Pi\,=\,\Pi^2 \quad\quad \Pi\,=\,\Pi^{\top}
\end{equation}
where $\Sigma^{1/2}$ is any square root of $\Sigma$ and $\Pi $ is an orthogonal projection matrix of rank $n$.

\begin{thm} \label{thmApprox}
The solutions of the   approximation problem  \eqref{Approx} are of the form 
$$
M=W\,W^{\top},  \qquad W = U_{n}  Q_{n}
$$
where $U_n$ is a  $N \times n$ matrix whose  columns are the first   $n$ normalized eigenvectors of  $\Sigma$, ordered  in the descending magnitude ordering of the corresponding first $n$ eigenvalues, collected in the diagonal matrix $\Lambda_n$, and  $Q_n$ is an arbitrary 
  $n\times n$ orthogonal matrix.
\end{thm}
\nd {\em Proof:}
Let
$ \Lambda := \diag \{\lambda_{1},\ldots,\lambda_{N}\}$ and $\Sigma\,=\, U\Lambda U^{\top}$ be the spectral decomposition of 
$\Sigma$ in which $U$ is an orthogonal matrix of eigenvectors. We can pick as a square rot of  $\Sigma$ the matrix $\Sigma^{1/2}\,:= \,U \Lambda^{1/2}$. 

Now, no matter how $\Sigma^{1/2}$ is chosen, the random vector 
${\mathbf e}:= \Sigma^{-1/2}\,{y}$ has orthonormal components. Hence using  \eqref{structM} the cost function of our minimum problem can be rewritten as 
\begin{align*}\label{lcl}
E\{\|{y} -M\, {y}\|^2\}\,&= \,\E\{\|\Sigma^{1/2}{\mathbf e} - \Sigma^{1/2}\, \Pi\,\Sigma^{-1/2} {\mathbf
y}\|^2\}\\
 &= \,\E\{\|\Sigma^{1/2} ({\mathbf e} -   \Pi \, {\mathbf e}\,)\|^2\} \\
&=\,\E\{\|\Lambda^{1/2} ({\mathbf e} -   \Pi \, {\mathbf e}\,)\|^2\} \\
 &=   \,\E\, ({\mathbf e} - \Pi\,{\mathbf e}\,)^{\top}
\Lambda \, ({\mathbf e} - \Pi\,{\mathbf e}\,)\\
&=\, \tr \left[ \Lambda\,\E ({\mathbf e} - \Pi\,{\mathbf e}\,)({\mathbf e} -
\Pi\,{\mathbf e}\,)^{\top}\right]
\end{align*}

where $\tr A := \sum a_{kk}$ is the trace of  $A$.  Our minimum problem can therefore be rewritten as
$$
\min_{\rank (\,\Pi\,)\, =\, n}\, \tr \{ \Lambda \Pi^{\perp}\}
$$
where $  \Pi^{\perp} := I - \Pi$  is the orthogonal projection onto the orthogonal complement of the subspace   $\im \Pi$.
 
Since the eigenvalues are ordered in decreasing order;  i.e. $\{\lambda_{1}\geq \ldots \geq \lambda_{t}\}$, one sees that the minimum of this function of $\Pi$ is reached when   $\Pi$ projects onto the subspace spanned by the first $n$ coordinate axes. In other words, $\Pi _{optimal} = \diag\{I_n ,\,0\}$ the minimum being 
$\lambda_{n+1}+\ldots +\lambda_{N}$. It is then evident that 
$$
M =\, U\Lambda^{1/2} \,\Pi _{optimal} \Lambda^{-1/2} U^{\top} = U_{n} U_{n}^{\top}
$$
Naturally, multiplying  $U_{n}$ by any orthogonal  $n \times n$ matrix does not change the result. \hfill $\Box$
 \medskip

Observe that  $\hat{y}=U_{n} U_{n}^{\top}y   $ is just the first $n$ {\em Principal Components Approximation}  of $y$.
In fact it is well-known that the PCA vector $\hat{y}$ can be seen as a linear transformation acting on  ${y}$ \cite{Hotelling-36}.
This result confirms in particular that the truncated PCA expansion   is optimal in the sense that it provides the best $M$ and the best approximation subspace for the criterion \eqref{Approx}. This characterization has been also found by a different technique studying  subspace approximation problems; see e.g. \cite{Yang-95}.

Note for future reference that the variance matrix of  $\hat{y}$ has rank $n$ since
\begin{equation} \label{Sigman}
\E\hat{y}\hat{y}^{\top}\!= \!U_{n} U_{n}^{\top}\! \E yy^{\top}U_{n} U_{n}^{\top}\!=\!U_{n} \diag \{\lambda_{1}, \ldots ,\lambda_{n} \}  U_{n}^{\top}
\end{equation}
and that this expression holds for arbitrary $N\geq n$. Naturally one should keep in mind that the eigenvector matrices now depend on $N$ but the eigenvalues stay fixed.

\section{Extension to infinite dimension}
In the same spirit,  consider now a stationary zero-mean process $\yb$ and any jointly stationary zero-mean process (both written as a doubly infinite column vectors) $\zb$, spanning a subspace $\Hb(\zb)\subset  \Hb(\yb)$ of dimension $n$.
 Any such process $\zb$ must be purely deterministic of rank $n$ and is then  uniquely determined  by any finite string of random variables $\{\zb(t)\}_{t\in I}$ induced on an interval $I$ of length $N\geq n$. This statement from \cite[page 276-277]{LPBook} is reported  for completeness in the following lemma.
 
 \begin{lemma}\label{Extension}
Any p.d. process $\zb$ of rank $n$  can be represented for all $t\in \Zbb$ by a state-space model (i.e. a stochastic realization) of the form
\begin{align}
	\xib(t+1)& = A \xib(t)   \label{stateeq} \\
	\zb(t) & = c\,\xib(t)   \label{outpeq}
	\end{align}
	where $\xib(t)= [\,\xi_1(t),\,\xi_2(t),\,\ldots, \xi_n(t)\,]^{\top}$ is an $n$-dimensional basis vector spanning the Hilbert space $\Hb(\zb_N)$ linearly generated by the $N \geq n$ random variables of the set $\{\zb(s)\,:\, t\geq s\geq t-N+1\}$, $A$  is a $n \times n$ unitary matrix and $c$ is a $n$-vector such that the pair $(A,c)$ is observable.
\end{lemma}
 \nd  {\em Proof:} See  \cite[p. 277]{LPBook}. \hfill $\Box$
 
 This linear system extends   in time the finite family of random variables $\{\zb(s)\}$, generators of $\Hb(\zb_N)$, to generate   the stationary p.d.  process $\zb$ defined on the whole time line $\Zbb$. Since this realization is uniquely determined by  the space $\Hb(\zb_N)$ modulo a choice of basis, it follows that the process $\zb$ is also uniquely determined by the space $\Hb(\zb_N)$. Hence all   covariances $\sigma(\tau)= \E\zb(t+\tau)\zb(t)$ are also uniquely defined and  determine the entries of the covariance operator  of the process. 
 \medskip

 In analogy to Problem \ref{ApproxProb}  let us  ask if there is a stationary process  $\zb$ spanning a subspace $\Hb(\zb)\subset \Hb(\yb)$ of dimension $n$,  which minimizes an approximation criterion of the type \eqref{Approx}.   If such a process exists we shall call it a {\em $n$-PC approximation of} $\yb$ and denote it by the symbol  $\hat\yb_n$.

Let then $I=[\, t, t+1, \ldots, t+N-1\,]$ denote a finite subinterval of the time line $\Zbb$ of length $N\geq n$ and  $X_I$ a $N\times \infty$   matrix with  entries $x_{j,k}$ equal to one if $k\in I$ and zero otherwise. Consider the finite random vectors $	\yb_N:=X_I\yb$ and $\zb_N:=X_I\zb$  and let $\hat\yb_N$  be  the random vector minimizing    the norm  $\E\{\|X_I\yb-X_I\zb\|^2\}$ for an   interval of length $N\geq n$. This problem is analogous to the problem \eqref{ApproxProb} where now $y$ is substituted by $\yb_N$. The solution  string determines by stationary extension (Lemma \ref{Extension})  a purely deterministic process $\hat\yb$ such that $\hat\yb_N:=X_I\hat\yb$ and   $\Hb(\hat\yb)= \Hb(\hat\yb_N)$ has dimension $n$.  Since $\Hb(\hat\yb)\subset  \Hb(\yb_N)\subset  \Hb(\yb)$,  this processes satisfies our requirements.   By stationarity   $\hat\yb$ is invariant with respect to translations of the interval $I$ provided its length $N$  is fixed. Then    $\hat\yb$ is a  $n$-PC approximation of $\yb$. The question now is to understand in what sense this approximation can get close to $\yb$ as $n\to \infty$.

 Since we are now  studying the behavior of  the $n$-PC approximation  of $\yb$ when  the dimension $n$ varies, we shall   attach a superscript  to $\hat\yb$  and  denote it by $\hat\yb^n$; likewise we shall do  to its covariance matrix, which will be denoted  $\hat\Sigmab^{n}$. 
 
 \begin{thm}\label{ProcessApprox}
 For each $n$ the $n$-PC approximation of $\yb$   has an (infinite) covariance operator $ \hat\Sigmab^n$ of rank $n$. The sequence  $\{ \hat\Sigmab^n\}$ converges weakly to $\Sigmab$ as $ n$  diverges to $\infty$, that is 
$$
 \psi^{\top} [\Sigmab -  \hat\Sigmab^n]\psi \rightarrow 0\qq \text{as}\q n\to \infty
 $$
for all functions  (row sequence) $\psi^{\top} := a^{\top}X_I$ having  support in $I$ where   $a \in \Rbb^N$ is arbitrary.
\end{thm}
\nd  {\em Proof:}
Consider  a  $n$-PC approximation $\hat\yb$ of  $\yb$ and the restriction  $X_I\hat\yb$ to  any  interval $I$ of length $N\geq n$. Recall that  $\hat\yb$ is now   denoted $\hat\yb^n$ and likewise for  its covariance matrix denoted  $\hat\Sigma^{n}$. By analogy to    \eqref{Sigman} the $N\times N$ truncation of this matrix has the structure 
\begin{equation} \label{SigmaN}
\hat\Sigma^{n}_N= U^{n}_N \diag \{\lambda_{1}, \ldots ,\lambda_{n} \} (U^{n}_N)^{\top}
\end{equation}
where the eigenvalues $\lambda_k$ are fixed and equal to the first  $n$ eigenvalues of  the $N\times N$-truncation of the covariance operator $\Sigmab$ of $\yb$. The $N\times n$ eigenvector matrices $U^{n}_N $ depend on $N$ as their dimension trivially grows with $N$.

We shall now show that all  $\hat\Sigma^{n}_N$ are the $N$-truncation of an {\em infinite} Toeplitz matrix $\hat\Sigmab^n$ of rank $n$ which is the covariance of the purely deterministic process $\hat\yb^n$. That this is so follows from the fact that  all  covariances $\hat\sigma^n(\tau)= \E\hat\yb^n(t+\tau)\hat\yb^n(t)$ are uniquely defined by the state space model \eqref{stateeq}, \eqref{outpeq} and  constitute  exactly  the entries of the (infinite) covariance operator $\hat\Sigmab^n$ of the p.d. process $\hat\yb^n$. Then, in particular, each finite matrix $\hat\Sigma^{n}_N$ must be a $N$-truncation of the same  $\hat\Sigmab^n$. Then from the expansion \eqref{SigmaN} it follows that this truncation  is just the (symmetric) $N\times N$  SVD approximation of rank $n$ of the  $N\times N$ truncation of $\Sigmab$, which is of course well defined for all finite $N$. 

 By a well-known property of the SVD (see e.g \cite[Chap. 2]{Golub-vL-89}) the variance matrix $\hat\Sigma^{n}_N$ of   $X_I\hat\yb^n$ is the symmetric $N\times N$ matrix of  rank $n$  which  has minimum distance from that of $X_I\yb$ in the Frobenius norm. This in turn implies that  the variance  $\hat\Sigmab^n$ is the (infinte) covariance matrix of rank $n$ which minimizes  
\begin{equation}\label{diff}
 \psi^{\top} (\Sigmab -  \hat\Sigmab^n)\psi =\psi^{\top} (\Sigmab_N -  \Sigma^{n}_N)\psi   .
 \end{equation}
for all functions $\psi$ having support in an interval $I\subset \Zbb$ of length $N\geq n$.
But the sequence \eqref{diff} is monotonically decreasing and nonnegative  for all $n$ and fixed  $N\geq n$. In fact by obvious properties of the SVD, for each $\psi$ of finite support of length $N$ it tends to zero with $n$ in a finite number of steps since when $n=N$ the difference is zero. But this happens for all $N$ and hence for all $\psi$ of finite support.\hfill $\Box$

\begin{rem}
Contrary to all submatrices $\Sigma_N$, the infinite covariance operator  $\Sigmab$ may not have eigenvalues (nor corresponding eigenvectors) and consequently the  
idea of PC approximation does not apply to the full matrix. For this reason the approximation  and the convergence results may not hold in a strong sense.
\end{rem}

\section{Approximation in the spectral domain}
We now go back to the problem formulation of Sec. \ref{Sec: ProbStat}. From \cite[Chap. 3]{LPBook} the processes $\yb$ and $\hat\yb^n$ have a spectral representation with random spectral measures $d Z(e^{i\theta})$ and $d Z^n(e^{i\theta})$ such that
\begin{align*}
\E dZ(e^{i\theta})dZ((e^{i\theta})^{*} &= \varphi(e^{i\theta})\Frac{d\omega}{2\pi} \,, \\
 \E dZ^n(e^{i\theta})dZ^n((e^{i\theta})^{*}& = \varphi_n(e^{i\theta})\Frac{d\omega}{2\pi}\,.
\end{align*}
Letting $\hat \psi(e^{i\omega}):=  \sum_{k=0}^{N-1} \psi(k) e^{i\omega k}$ one has 
\begin{align*}
 \psi^{\top} \Sigmab \psi &= \E[ \sum_{k=0}^{N-1} \psi(k) \yb(k)\,]^2= \E [ \int_{-\pi}^{\pi} \hat \psi(e^{i\omega}) dZ(e^{i\theta})\,]^2 \\
 &= \int_{-\pi}^{\pi} \hat \psi(e^{i\omega}) \varphi(e^{i\omega})\hat \psi(e^{i\omega})^{*}\Frac{d\omega}{2\pi}
\end{align*}
and similarly for $ \psi^{\top} \Sigmab^n \psi$ which can be written
\begin{align*}
 \psi^{\top} \Sigmab^n \psi &= \E[ \sum_{k=0}^{N-1} \psi(k) \hat\yb(k)\,]^2\\
 &= \E [ \int_{-\pi}^{\pi} \hat \psi(e^{i\omega}) dZ^n(e^{i\theta})\,]^2\\
 & = \int_{-\pi}^{\pi} \hat \psi(e^{i\omega}) \varphi_n(e^{i\omega})\hat \psi(e^{i\omega})^{*}\Frac{d\omega}{2\pi}
\end{align*}
Therefore \eqref{Wconv} follows from Theorem \ref{ProcessApprox}.

Now note that, because of the orthogonality \eqref{OrthM}, $\psi^{\top} (\yb - \hat\yb^n) \perp  \psi^{\top}\hat\yb^n$ and  the difference \eqref{diff} can be rewritten  as  
$ 
\E [\psi^{\top} (\yb - \hat\yb^n)]^2
$ 
which   also must tend to  zero when $n\to \infty$ for all functions $\psi$ having support in an interval $I\subset \Zbb$ of length $N\geq n$. Therefore $\hat\yb^n$ converges weakly to $\yb$.

\bibliographystyle{IEEEtran}
\bibliography{biblio-FreqEst,Toeplitz,Bookrefs}

\end{document}